\documentclass[12pt]{article}
\usepackage{amsmath,amssymb,amsthm,amscd}
\usepackage[all]{xy}

\newcommand{\Z}{\mathbb{Z}}
\newcommand{\Q}{\mathbb{Q}}

\newcommand{\T}{\mathbb{T}}

\newtheorem{defn}{Definition}[section]
\newtheorem{thm}[defn]{Theorem}
\newtheorem{ex}[defn]{Example}

\newtheorem{cor}[defn]{Corollary}

\newtheorem{rem}[defn]{Remark}

\title{A homology theory for Smale spaces: a summary}
\author{Ian F. Putnam\thanks{Supported in part by a
grant from NSERC, Canada},\\
Department of Mathematics and Statistics,\\
University of Victoria,\\
Victoria, B.C., Canada V8W 3P4}
\date{ }

\begin{document}

\maketitle

\begin{abstract}
We consider Smale spaces, a particular class of hyperbolic topological dynamical systems, which include the basic sets for Smale's Axiom A systems. We present a homology theory for such systems which is based on Krieger's dimension group in the special case of shifts of finite type. This theory provides a Lefschetz formula relating trace data with the number of periodic points of the system.
\end{abstract}

\section{Introduction}

Smale introduced the notion of an Axiom A diffeomorphism of a compact manifold \cite{S:diffdyn}. 
For such a system, a basic set
is an invariant subset of the non-wandering set which is irreducible in a certain sense. One of
Smale's key observations was that  such a set need not be a submanifold. Typically, it is some type of fractal object.
The aim of this article is to introduce a new type of homology theory for such spaces. This takes as
its starting point the notion of the dimension group of a shift of finite type introduced by
Krieger \cite{Kr:dimgp} and the fundamental result of Bowen \cite{Bow:markpart} that every basic 
set is a factor  of a shift of finite
type.

In an effort to give a purely topological (i.e. without reference to any smooth structure) description of the dynamics on a basic set,  Ruelle introduced the notion 
of a Smale space \cite{R:ThermForm}: a Smale space is a compact metric space, $(X,d)$, and a homeomorphism,
$\varphi$, of $X$, which possesses canonical coordinates of contracting and expanding directions. 
The precise definition involves the existence of a map $[,]$ giving canonical coordinates. Here, we review only the features necessary for the statements of our results.

There is a constant $\epsilon_{X} > 0$ and, for each $x$ in $X$ and $0 < \epsilon \leq \epsilon_{X}$, there  are sets $X^{s}(x, \epsilon)$ and $X^{u}(x, \epsilon)$, called the local stable and unstable sets, respectively, whose product is homeomorphic to a neighbourhood of $x$. As $\epsilon$ varies, these form a neighbourhood base at $x$. Moreover, there is  a constant $0 < \lambda < 1$ such that 
\begin{eqnarray*}
 d(\varphi(y),\varphi(z)) & \leq \lambda d(y,z), & y, z \in X^{s}(x, \epsilon_{X})\\
d(\varphi^{-1}(y),\varphi^{-1}(z)) & \leq \lambda d(y,z), & y, z \in X^{u}(x, \epsilon_{X})
\end{eqnarray*}
The bracket $[x, y]$ is the unique point in the intersection of 
$X^{s}(x, \epsilon_{X})$ and $X^{u}(y, \epsilon_{X})$.
We say that $(X, \varphi)$ is non-wandering if every point of $X$ is non-wandering for $\varphi$ \cite{KH:book}.

Stable and unstable equivalence relations are defined by
\begin{eqnarray*}
 R^{s} & = & \{ (x,y) \mid \lim_{n \rightarrow + \infty } d(\varphi^{n}(x), \varphi^{n}(y)) = 0 \} \\
R^{u} & = & \{ (x,y) \mid \lim_{n \rightarrow + \infty } d(\varphi^{-n}(x), \varphi^{-n}(y)) = 0 \}.
\end{eqnarray*}
We let $X^{s}(x)$ and $X^{u}(x)$ denote the stable and unstable equivalence classes of  $x$ in $X$.

The main examples of such systems are shifts of finite type (of which we will say more in a moment), hyperbolic toral automorphisms, solenoids, substitution tiling spaces (under some hypotheses) and, most importantly, the basic sets for Smale's Axiom A systems \cite{S:diffdyn,Bow:CBMS}.

Let $(Y, \psi)$ and  $(X, \varphi)$ be  Smale spaces. A 
\emph{factor map} from $(Y, \psi)$ to $(X, \varphi)$ is a function $\pi:Y \rightarrow X$ which is continuous, surjective and satisfies $\pi \circ \psi = \varphi \circ \pi$. It is clear that, for any $y$ in $Y$, $\pi(X^{s}(y)) \subset X^{s}(\pi(y))$ and $\pi(X^{u}(y)) \subset X^{u}(\pi(y))$.
David Fried \cite{Fr:finpres} defined $\pi$ to be $s$-resolving (or $u$-resolving)
if, for every $y$ in $Y$, the restriction of $\pi$ to $X^{s}(y)$ (or to $X^{u}(y)$, respectively) is injective.  We say that
$\pi$ is \emph{$s$-bijective} (or \emph{$u$-bijective}) if, for every $y$ in $Y$, $\pi$ is a bijection from $X^{s}(y)$ to $X^{s}(\pi(y))$ (or from $X^{u}(y)$ to $X^{u}(\pi(y))$, respectively). This actually implies that $\pi$ is a local homeomorphism from the local stable sets (or unstable sets, respectively) in $Y$  to those in $X$. In the case that $(X, \varphi)$ is  non-wandering, $s$-resolving ($u$-resolving) and $s$-bijective ($u$-bijective, respectively) are equivalent.

\section{Shifts of Finite type}

Shifts of finite type are described in detail in \cite{LM:SymDynbook}. 
We consider a finite directed graph $G$. This consists of a finite vertex set $G^{0}$, a finite edge set $G^{1}$ and maps $i, t$ (for initial and terminal) from $G^{1}$ to $G^{0}$.  The associated shift space
\[
 \Sigma_{G} = \{ e = (e^{k})_{k \in \Z} \mid e^{k} \in G^{1}, t(e^{k}) = i(e^{k+1}), k \in \Z \}
\]
consists of all bi-infinite paths in $G$, specified as an edge list. The map $\sigma$ is the left shift on $\Sigma_{G}$ defined by $\sigma(e)^{k} = e^{k+1}$, for all $e$ in $\Sigma_{G}$ and $k$ in $\Z$.
 By a shift of finite type, we mean any system topologically conjugate to $(\Sigma_{G}, \sigma)$, for some graph $G$. This is not the usual definition, but is equivalent to it (see Theorem 2.3.2 of \cite{LM:SymDynbook}).
We observe that such systems are Smale spaces by noting that the bracket 
operation is defined as follows. For $e,f$ in $\Sigma_{G}$, $[e,f]$ is defined if $t(e^{0}) = t(f^{0})$ and then it is the sequence,  $( \ldots, f^{-1}, f^{0}, e^{1}, e^{2}, \ldots)$. For any $l \geq 1$ and $e_{0}$ in $\Sigma_{G}$, the local stable and unstable sets are given by 
\begin{eqnarray*}
\Sigma_{G}^{s}(e_{0}, 2^{-l}) & =  & \{ e \in \Sigma \mid e^{k} = e_{0}^{k}, k \geq -l \}, \\
 \Sigma_{G}^{u}(e_{0}, 2^{-l}) &  =  & \{ e \in \Sigma \mid e^{k} = e_{0}^{k}, k \leq l \}.
\end{eqnarray*}
\begin{thm} Shifts of finite type are exactly the zero-dimensional (i.e. totally disconnected) Smale spaces.
\end{thm}

The fundamental r\^{o}le of shifts of finite type is demonstrated by the following universal property, due to  Bowen \cite{Bow:markpart}.

\begin{thm}[Bowen]
\label{thm:Bowen}
 Let $(X, \varphi)$ be a non-wandering Smale space. Then there exists a non-wandering shift of finite type
$(\Sigma, \sigma)$ and a finite-to-one factor map
\[
 \pi: (\Sigma, \sigma) \rightarrow (X, \varphi).
\]
\end{thm}

\section{Krieger's dimension group invariant}

Krieger \cite{Kr:dimgp}  defined the past and future dimension groups of a shift 
 of finite, $(\Sigma, \sigma)$, as follows. Consider 
$\mathcal{D}^{s}(\Sigma, \sigma)$ to be the collection
of compact open subsets of $\Sigma^{s}(e, \epsilon)$, as $e$ varies
over $\Sigma$ and $0 < \epsilon \leq \epsilon_{X}$.
We let $\sim$ denote the smallest equivalence relation on $\mathcal{D}^{s}(\Sigma,\sigma)$
such that 
\begin{enumerate}
 \item If $E, F$ are in  $\mathcal{D}^{s}(\Sigma, \sigma)$ with
$[E, F] = F, [F, E] = E$ (meaning both are defined),
then $E \sim F$,
\item If $E, F, \varphi(E)$ and  $\varphi(F)$ are all
in $\mathcal{D}^{s}(\Sigma, \sigma)$, then
$E \sim F$ if and only if $\varphi(E) \sim \varphi(F)$.
\end{enumerate}
We generate a free abelian group on the equivalence classes 
of elements, $E$, of $\mathcal{D}^{s}(\Sigma, \sigma)$
 (denoted $[E]$) subject to the additional 
relation that $[E \cup F ] = [E] +[F]$, if $E \cup F$ is in 
 $\mathcal{D}^{s}(\Sigma, \sigma)$ and $E$ and $F$ are disjoint.
The result is denoted by $D^{s}(\Sigma, \sigma)$. There is an 
analogous definition of $D^{u}(\Sigma, \sigma)$.

\begin{thm}[Krieger \cite{Kr:dimgp}]
 Let $G$ be a finite directed graph and $(\Sigma_{G}, \sigma)$ be the 
associated shift of finite type. Then $D^{s}(\Sigma_{G}, \sigma)$ 
(or $D^{u}(\Sigma_{G}, \sigma)$, respectively) is 
isomorphic to the inductive limit of the sequence
\[
 \Z G^{0} \stackrel{\gamma^{s}}{\rightarrow} \Z G^{0} \cdots 
\]
where $\Z G^{0}$ denotes the free abelian group on the vertex set
$G^{0}$ and the map $\gamma(s) = \sum_{t(e)=v} i(e)$, for
any $v$ in $G^{0}$ ( or respectively, replacing $\gamma^{s}$ by
$\gamma^{u}$, whose definition is the same, interchanging the r\^{o}les
of $i$ and $t$).
\end{thm}

The first crucial result for the development of our  theory
is the following  functorial property of $D^{s}$ and $D^{u}$, which can be found in 
\cite{BMT:Memoir}.

\begin{thm}
\label{thm:D_functor}
 Let $(\Sigma, \sigma)$ and $(\Sigma', \sigma)$ be shifts of finite
type and let \newline
$\pi: (\Sigma, \sigma) \rightarrow (\Sigma', \sigma)$
be a factor map.
\begin{enumerate}
 \item If $\pi$ is  $s$-bijective, then there are
natural homomorphisms
\[
 \begin{array}{ccccc}
 \pi^{s} & : & D^{s}(\Sigma, \sigma) & \rightarrow & 
                        D^{s}(\Sigma', \sigma), \\
\pi^{u*} & : & D^{u}(\Sigma', \sigma) & \rightarrow & 
                        D^{u}(\Sigma, \sigma).
 \end{array}
\]
 \item If $\pi$ is $u$-bijective, then there are
natural homomorphisms
\[
 \begin{array}{ccccc}
 \pi^{u} & : & D^{u}(\Sigma, \sigma) & \rightarrow & 
                        D^{u}(\Sigma', \sigma), \\
\pi^{s*} & : & D^{s}(\Sigma', \sigma) & \rightarrow & 
                        D^{s}(\Sigma, \sigma).
 \end{array}
\]
\end{enumerate}
\end{thm}

The idea is simple enough: in the covariant case, the induced 
map sends the class of a set $E$ in $\mathcal{D}^{s}(\Sigma, \sigma)$
to the class of $\pi(E)$, while in the contravariant case,  the map sends the
class of  $E'$ in $\mathcal{D}^{s}(\Sigma', \sigma)$
to the class of  $\pi^{-1}(E')$. The latter is not correct since
$\pi^{-1}(E')$ may not even be contained in a single
stable equivalence class, but it suffices that it may
be written as a finite union of elements of 
$\mathcal{D}^{s}(\Sigma, \sigma)$. These ideas work correctly under the stated  hypotheses.

\section{$s/u$-bijective pairs}

The key ingredient in our construction is the following notion.

\begin{defn}
 Let $(X, \varphi)$ be a Smale space. An \emph{$s/u$-bijective pair}, $\pi$,
for $(X, \varphi)$ consists of Smale spaces $(Y, \phi)$ and 
$(Z, \zeta)$ and factor maps 
\[
\pi_{s}: (Y, \psi) \rightarrow (X, \varphi), \hspace{.5cm}
\pi_{u}: (Z, \zeta) \rightarrow (X, \varphi)
\]
such that
\begin{enumerate}
 \item $Y^{u}(y, \epsilon)$ is totally disconnected,
for all $y$ in $Y$ and $ 0 < \epsilon \leq \epsilon_{Y}$,
\item $\pi_{s}$ is $s$-bijective,
\item $Z^{s}(z, \epsilon)$ is totally disconnected,
for all $z$ in $Z$ and $ 0 < \epsilon \leq \epsilon_{Z}$,
\item $\pi_{u}$ $u$-bijective.
\end{enumerate}
\end{defn}

To summarize the idea in an informal way, the space
$Y$ is an extension of $X$, where the local unstable sets
are totally disconnected, while the local stable sets are
homeomorphic to those in $X$. The existence of such 
$s/u$-bijective pairs, at least for non-wandering $(X, \varphi)$, can be 
deduced from the results of \cite{P:lifting}. It
 can be viewed as a coordinate-wise version of Bowen's theorem.

\begin{thm}
If $(X, \varphi)$ is non-wandering, then there exists 
an $s/u$-bijective  pair for $(X, \varphi)$.
\end{thm}

\begin{defn}
 Let $\pi = (Y, \psi, \pi_{u}, Z, \zeta, \pi_{s})$ be an $s/u$-bijective pair
for \newline
$(X, \varphi)$. For each $L, M \geq 0$, we define
\begin{eqnarray*}
 \Sigma_{L,M}(\pi) & =  & \{ (y_{0}, \ldots, y_{L}, z_{0}, \ldots , z_{M})  \mid      y_{l} \in Y, z_{m} \in Z, \\
   &   &   \pi_{s}(y_{l}) = \pi_{u}(z_{m}),
       0 \leq l \leq L, 0 \leq m \leq M \}
\end{eqnarray*}
For simplicity, we also denote $\Sigma_{0,0}(\pi)$ by $\Sigma(\pi)$.

We define $\sigma_{L,M}: \Sigma_{L,M}(\pi) \rightarrow \Sigma_{L,M}(\pi)$
by 
\[
\sigma_{L,M}(y_{0}, \ldots, y_{L}, z_{0}, \ldots, z_{M})
  =   (\psi(y_{0}), \ldots, \psi(y_{L}), \zeta(z_{0}), \ldots, 
     \zeta(z_{M})).
\]

For $L \geq 1$ and $0 \leq l \leq L$, we let $\delta_{l,}: \Sigma_{L,M}(\pi) \rightarrow
\Sigma_{L-1,M}(\pi)$ be the map which deletes entry $y_{l}$. Similarly,
the map $\delta_{,m}: \Sigma_{L,M}(\pi) \rightarrow \Sigma_{L,M-1}(\pi)$ deletes entry $z_{m}$, for $M \geq 1, 
0 \leq m \leq M$.
\end{defn}

The important properities of these systems and maps is summarized as follows. 

\begin{thm}
For every $L,M \geq 0$, $(\Sigma_{L,M}(\pi), \sigma_{L,M})$ is a shift of finite type.
 For $L \geq 1$, $0 \leq l \leq L$, the map $\delta_{l,}$ is an
 $s$-bijective factor map. For $M \geq 1$, $0 \leq m \leq M$, the map $\delta_{,m}$ is a
$u$-bijective factor map.
\end{thm}

\section{Homology}

There are actually two homology theories here. One, based on the dimension group $D^{s}$ will be denoted by $H^{s}_{*}$ and the other, based on $D^{u}$, will be denoted by $H^{u}_{*}$. We will concentrate on the former for the remainder of this note.

\begin{defn}
Let $(X, \varphi)$ be  a Smale space and suppose that \newline
$\pi = (Y, \psi, \pi_{s}, Z, \zeta, \pi_{u})$ is an $s/u$-bijective pair
for $(X, \varphi)$.
We define
\[
 C_{N}^{s}(\pi) = \oplus_{L-M= N} D^{s}(\Sigma_{L,M}(\pi), \sigma_{L,M})
\]
for every $N$ in $\Z$ and a boundary map $\partial^{s}_{N}(\pi): C_{N}^{s}(\pi) \rightarrow C_{N-1}^{s}(\pi)$ by
\[
 \partial^{s}_{N}(\pi) | D^{s}(\Sigma_{L,M}, \sigma_{L,M})
 = \sum_{l = 0}^{L} (-1)^{l}
   \delta_{l,}^{s} + \sum_{m=0}^{M+1} (-1)^{m + L} 
   \delta_{,m}^{s*}.
\]
where, in the special case $L = 0$, we set 
$\delta_{0,}^{s} = 0$. 

We define $ H^{s}_{*}(\pi) $ to be the homology of the complex
$(C_{*}^{s}(\pi), \partial_{*}^{s}(\pi))$.
\end{defn}

We want to establish some basic properties of 
our theory.
The first crucial result   is the following. It is stated in a slightly informal manner, but it conveys the main idea.

\begin{thm}
 $H^{s}_{*}(\pi)$ is independent of the $s/u$-bijective pair \newline 
$\pi = (Y, \pi_{s}, Z, \pi_{u})$ and depends only on $(X, \varphi)$. 
\end{thm}

In light of this, we denote $H^{s}_{*}(\pi)$ by $H^{*}_{N}(X, \varphi)$ instead. It is defined provided that there exists 
an $s/u$-bijective pair for $(X, \varphi)$, which is true for all  non-wandering Smale spaces.

\begin{thm}
\label{thm:H_funct}
The homology theory $H^{s}_{*}$  is functorial in the following sense: if $\pi: (Y, \psi) \rightarrow (X, \varphi)$ is an $s$-bijective factor map, then there are induced group homomorphisms
\[
 \pi^{s}: H^{s}_{N}(Y, \psi) \rightarrow H^{s}_{N}(X, \varphi), 
\]
for all $N$ in $\Z$. If the map $\pi$ is a $u$-bijective factor, then there are induced group homomorphisms
\[
 \pi^{s*}:  H^{s}_{N}(X, \varphi)  \rightarrow H^{s}_{N}(Y, \psi), 
\]
for all $N$ in $\Z$.
\end{thm}

To describe the next property,  we assume that we have chosen a graph
$G$ which presents the shift of finite type $(\Sigma(\pi), \sigma)$ associated
with our $s/u$-bijective pair. More specifically, for each 
$(y, z)$ in $\Sigma(\pi)$ and each integer $k$, $e^{k}(y,z)$
is an edge in $G$ and the map sending $(y, z)$ in $\Sigma(\pi)$
to $(e^{k}(y, z))_{k \in \Z}$ is a conjugacy between 
$(\Sigma(\pi), \sigma)$ and $(\Sigma_{G}, \sigma)$. We make the additional assumption
that the presentation
 is \emph{regular} in the following sense. For any $(y,z), (y',z')$ in $\Sigma(\pi)$ such that
$t(e^{0}(y,z)) = t(e^{0}(y',z'))$, it follows that
$[\pi_{s}(y), \pi_{s}(y')]$ is defined in $X$,
$[y,y']$ is defined in $Y$ and $[z, z']$ is defined in $Z$  and 
we have 
\[
 \pi_{s}[y, y'] = \pi_{u}[z,z'] = [\pi_{s}(y), \pi_{s}(y')].
\]
This may always be achieved by passing to a higher block presentation
of $G$ \cite{LM:SymDynbook}. As a consequence, it is fairly easy to see that,
for all $L, M \geq 0$, there exists a graph $G_{L,M}$ consisting
of $L+1$ by $M+1$ arrays of vertices and edges from $G$ such that
the map sending $(y_{0}, \ldots, y_{L}, z_{0}, \ldots, z_{M})$
in $\Sigma_{L,M}(\pi)$ to the sequence $(e^{k}(y_{l}, z_{m}))_{k \in \Z}$ 
in a conjugacy with $\Sigma_{G_{L,M}}$.

We note that the groups $S_{L+1}$ and $S_{M+1}$ act on both
$\Sigma_{L, M}(\pi)$ and $G_{L,M}$.  We let $\Z [ G^{0}_{L,M})$ be the quotient 
$\Z G^{0}_{L,M}$ by the relations that $v = 0$ if two distinct rows of
$v$ are equal and $\alpha \cdot v = sgn(\alpha) v$, for any
$v$ in $G^{0}_{L,M}$ and $\alpha$ in $S_{L+1}$. The map $\gamma^{s}_{L,M}$
associated with $G_{L,M}$ descends to a well-defined map on the quotient and 
we let $D^{s}[G_{L,M})$ denote the limit group.
It must also be checked that  the boundary map descends to the complex
$C^{s}_{N}[\pi) = \oplus_{L-M=N} D^{s}[G_{L,M})$,
and we denote the result by $\partial^{s}[\pi)$.
There analogous complexes $(C^{s}(\pi], \partial^{s}(\pi])$ which treats the
actions of the various $S_{M+1}$'s in an analogous way. In this case, it is a subcomplex
of our original. Finally, there is a quotient of this,  denoted
$(C^{s}[\pi], \partial^{s}[\pi])$ which includes both $S_{L+1}$ and $S_{M+1}$ actions.

\begin{thm}
 The natural maps between the   complexes \newline
$(C^{s}(\pi), \partial^{s}(\pi))$, $(C^{s}[\pi), \partial^{s}[\pi))$,
$(C^{s}(\pi], \partial^{s}(\pi])$ and 
$(C^{s}[\pi], \partial^{s}[\pi])$ all induce isomorphisms at
the level of homology.
\end{thm}

This result is very useful for computational purposes
since the groups $\Z [ G^{0}_{L,M}]$ have many fewer
generators  $\Z G^{0}_{L,M}$. In addition,
every $s$-bijective or $u$-bijective map is finite-to-one, meaning that 
there is a uniform upper bound on the cardinality of a pre-image. Hence, there exists
some $N_{0}$ such that, if either $L$ or $M$ exceed $N_{0}$, any element of
$\Sigma_{L,M}(\pi)$ will contain a repeated entry. In consequence,
the group $\Z[G^{0}_{L,M} ]$, for all such $L, M$, will be trivial.

\begin{cor}
 For any Smale space $(X, \varphi)$ which has an $s/u$-bijective pair, the 
groups $H_{N}^{s}(X, \varphi)$ are finite rank and are non-zero
for only finite many values of $N$.
\end{cor}

One of the most important aspects of Krieger's theory is that, for 
any shift of finite type $(\Sigma, \sigma)$, the group
$D^{s}(\Sigma, \sigma)$ has a natural order structure. This carries over
to our theory as follows.

\begin{thm}
 If $(X, \varphi)$ is a non-wandering Smale space, then the group
\newline
$H^{s}_{0}(X, \varphi)$ has a natural order structure, which is preserved
under the functorial properties of Theorem \ref{thm:H_funct}.
\end{thm}

Finally, we have the following analogue of the Lefschetz formula. 
Given $(X, \varphi)$, we can regard $\varphi$ as a 
factor map from this system to itself. It is both  $s$-bijective 
and $u$-bijective and so, by Theorem \ref{thm:H_funct}, induces an automorphism of our invariant,
denoted $\varphi^{s}$. 
The following result, already known in the case of shifts of finite type, 
uses ideas of Manning \cite{Man:ratzeta}.

\begin{thm}
 For any non-wandering Smale space $(X, \varphi)$ and $p \geq 1$, we have
\begin{eqnarray*}
 & \sum_{N \in \Z} (-1)^{N} Tr [ (\varphi^{s})^{p}: H^{s}_{N}(X, \varphi) \otimes \Q \rightarrow H^{s}_{N}(X, \varphi) \otimes \Q ]   &  \\
    & =   \# \{ x \in X \mid 
\varphi^{p}(x) = x \}. & 
 \end{eqnarray*}
\end{thm}

\section{Examples}

We present four examples where the  computations above may be carried out quite explicitly. The full details 
of the last three are in preparation \cite{BP:egs}. All of the examples are computed 
using the  double complex $C^{s}[\pi]$; it is that one which is 
described in each example.

\begin{ex}
 Suppose $(\Sigma, \sigma)$ is a shift of finite type.
In this case, an $s/u$-bijective pair is just $(Y, \psi) = (Z, \zeta) = (\Sigma, \sigma)$. Only the $0,0$-term in the double complex is non-zero and it is just $D^{s}(\Sigma, \sigma)$. Hence, $H^{s}_{N}(\Sigma, \sigma)$ is just 
$D^{s}(\Sigma, \sigma)$, for $N = 0$, and zero otherwise.
\end{ex}

\begin{ex}
 For $m \geq 2$, let $(X, \varphi)$ be the $m^{\infty}$-solenoid.
More specifically, we let
\[
 X = \{ (z_{0}, z_{1}, \ldots ) \mid z_{n} \in \T, z_{n} = z_{n+1}^{m}, n \geq 0 \},
\]
with the map
\[
 \varphi(z_{0}, z_{1}, \ldots ) = (z_{0}^{m}, z_{1}^{m}, z_{2}^{m}, \ldots ),
\]
for $(z_{0}, z_{1}, \ldots )$ in $X$. 
In this case, there is an $s$-bijective factor map onto $(X, \varphi)$ 
from the full $m$-shift (i.e. $G$ is the graph with one vertex and 
$m$  edges, although it is necessary to pass to a higher block presentation for the map to be regular). The 
$s/u$-bijective
 pair here is $(Y, \psi) = (\Sigma_{G}, \sigma)$ and $(Z, \zeta) = (X, \varphi)$. The only non-zero groups in the
 double complex occur for $L = 0 = M$ and $L=1, M=0$ and these are $\Z[m^{-1}]$ and $\Z$, respectively.
 The boundary maps are all zero (only one needs to be computed) and $H^{s}_{N}(X, \varphi)$ is isomorphic to $\Z[ m^{-1}]$, 
for $N=0$,  $\Z$, for $N=1$ and zero for all other $N$.
\end{ex}

\begin{ex}
 Let $X$ be the $2$-torus, $\T^{2}$, and $\varphi$ be the hyperbolic automorphism determined by the matrix
\[
A = \left( \begin{array}{cc} 1 & 1 \\ 1 & 0 \end{array} \right).
\]
In this case, a Markov partition can be chosen so that the natural quotient factors in two ways, as an $s$-bijective followed by $u$-bijective and vice verse. (In fact, it is the Markov partition with three rectangles which appears in many dynamics texts.) The only non-zero terms in the double complex are in positions $(L,M) = (0,0), (1,0), (0,1), (1,1)$. The calculation yields $H^{s}_{N}(X, \varphi)$ is $\Z$ for $N = 1$ and $N=-1$ and is $\Z^{2}$ for $N=0$. Notice that the homology coincides with that of the torus, except with a dimension shift.
\end{ex}

\begin{ex} There is an example, $(X, \varphi)$, roughly based on the Sierpinski gasket. We do not give any details except to mention that it is not a shift of finite type, but its homology is the same as the full $3$-shift. 
\end{ex}

\section{Concluding remarks}

\begin{rem}
 It is certainly a natural question to ask whether this theory can be computed from other (already existing) machinery. A more specific question would be to relate our homology to, say, the Cech cohomology of the classifying space of the topological equivalence relation $R^{s}$. (For a discussion of the topology, see \cite{{P:Smale_C*}}.)  There are examples, such as the first three above, where they are different, but only up to a dimension shift (depending on the space under consideration).
\end{rem}

\begin{rem}
 An important motivation in the construction of this theory was to compute the $K$-theory of certain $C^{*}$-algebras associated with the Smale space $(X, \varphi)$. See \cite{{P:Smale_C*}} for a discussion of these  $C^{*}$-algebras. At present, there seems to be a spectral sequence which relates the two; this work is still in progress.
\end{rem}

\bibliographystyle{amsplain}

\end{document}